\documentclass[a4paper,reqno,10pt]{amsart}
\usepackage{amssymb}
\usepackage{amsfonts}
\usepackage{graphicx}
\usepackage{color}
\usepackage{epstopdf}

\newtheorem{theorem}{Theorem}
\theoremstyle{plain}

\newtheorem{lemma}{Lemma}

\newtheorem{proposition}{Proposition}
\newtheorem{remark}{Remark}

\renewcommand\bigskip{\medskip}

\def\to{\rightarrow}
\def\bc{\begin{center}}
\def\ec{\end{center}}
\def\be{\begin{equation}}
\def\ee{\end{equation}}
\def\P{\mathbb P}

\def\N{\mathbb N}

\def\R{\mathbb R}

\def\E{\mathbb E}

\begin{document}
\title[]{Multifractal analysis of some multiple ergodic averages in linear Cookie-Cutter dynamical systems}

\author[A. H. Fan]{AiHua Fan}
\address{A. F. Fan: 
 LAMFA, UMR 7352  CNRS, Universit\'e de Picardie, 33 Rue
Saint Leu, 80039 Amiens, France. 
E-mail: ai-hua.fan@u-picardie.fr}

\author[L. M. Liao]{LingMin Liao}
\address{L. M. Liao: LAMA UMR 8050 CNRS, Universit\'e Paris-Est Cr\'eteil,
61 Avenue du G\'en\'eral de Gaulle, 94010 Cr\'eteil Cedex, France.
E-mail: lingmin.liao@u-pec.fr}

\author[M. Wu]{Meng Wu}
\address{M. Wu:  Department of Mathematical Sciences, P.O. Box 3000, 90014 University of Oulu,
Finland.
E-mail: meng.wu@oulu.fi}

\begin{abstract}
In this paper,  we study the multiple ergodic averages of a locally constant real-valued function
in linear Cookie-Cutter dynamical systems. The multifractal spectrum of these multiple ergodic averages is completely determined.
\end{abstract}

\maketitle

\markboth{Multifractal analysis of some multiple ergodic averages}{Ai-Hua FAN, Lingmin LIAO and Meng WU}%

\section{Introduction and statement of results }\label{int}

Let $T : X\to X$ be a continuous map on a compact metric space $X$. Let $f_{1},\cdots,f_{\ell}$  $(\ell\geq 2)$ be $\ell$ bounded real-valued 
functions on $X$. The following {\em multiple ergodic average}
$$\frac{1}{n}\sum_{k=1}^nf_{1}(T^kx)f_2(T^{2k}x)\cdots f_{\ell}(T^{\ell k}x)$$
is widely studied in ergodic theory by Furstenberg \cite{Furstenberg},   Bourgain \cite{Bourgain}, Host and Kra \cite{HK}, Bergelson, Host and Kra \cite{BHK} and others. 
Fan, Liao and Ma \cite{FLM} and Kifer \cite{Kifer}  have independently studied such multiple ergodic averages from the
point of view of multifractal analysis.
Later on,  the multifractal analysis of multiple ergodic averages have attracted much attention.  
First works are done on symbolic spaces.  
Let $m\geq 2$ be an integer and $S=\{0,\cdots,m-1\}$. Consider the symbolic space $\Sigma_m=S^{\N^*}$ endowed with the metric 
$$d(x,y)=m^{-\min\{n,x_n\neq y_n\}}, \ \ \forall x, y\in \Sigma_m.$$ 
The first object of study was the Hausdorff dimension of the following level sets (\cite{FLM})
$$E(\alpha)=\left\{(x_k)_{k=1}^\infty\in \Sigma_2\ : \ \lim_{n\to \infty}\frac{1}{n}\sum_{k=1}^n x_kx_{2k}=\alpha\right\},\ \ \alpha\in [0,1].$$
More generally we may consider the Hausdorff spectrum of the  following level sets of multiple ergodic averages
 \begin{equation}\label{equ 1}
 E^\ell_\varphi(\alpha)=\left\{(x_k)_{k=1}^\infty\in \Sigma_m\ : \ \lim_{n\to \infty}\frac{1}{n}\sum_{k=1}^n \varphi (x_k,x_{kq},\cdots, x_{kq^{\ell-1}})=\alpha\right\}, \alpha\in \R
 \end{equation}
 where $q\geq 2,\ell\geq 2$ are integers and $\varphi$ is a real-valued function defined on $\{0,\cdots,m-1\}^\ell$. The level set $E(\alpha)$ then corresponds to the set $E^\ell_\varphi(\alpha)$ with special choice $q=2,\ell=2$ and $\varphi(x,y)=xy$. 
See the works of Kenyon, Peres and Solomyak \cite{KPS,KPS1}, Peres, Schmeling, Seuret and Solomayk \cite{PSSS} on some specific subsets of level sets $E(\alpha)$.  
See Peres and Solomyak \cite{PS} for 
the multifractal analysis of $E(\alpha)$. 
Fan, Schmeling and Wu \cite{FSW,FSW1} have considered a class of functions $\varphi$ that are involved in (\ref{equ 1}). 
Fan, Schmeling and Wu \cite{FSW_V} have also considered some similar averages called $V$-statistics.

All of the above mentioned results concentrated on the full shift dynamical system  $(\Sigma_m, \sigma)$ where the Lyapunov exponent of the shift transformation is constant. Recently, Liao and Rams \cite{LR} performed the multifractal analysis of a class of special multiple ergodic averages for some systems with non-constant Lyapunov exponents. 
More precisely, they considered a piecewise linear map $T$ on the unit interval with two branches. Let $I_0,I_1\subset[0,1]$ be two intervals with disjoint interiors. Suppose that for each $i\in \{0,1\}$, the restriction $T: I_i\to [0,1]$ is bijective and linear with slop $e^{\lambda_i}$, $\lambda_i >0$. Let $J_T$ be the {\em repeller} of $T$, i.e.
$$J_T:=\bigcap_{n=1}^\infty T^{-n}[0,1].$$ 
Then $(J_T,T)$ becomes a dynamical system. As in \cite{FLM,FSW,PS}, Liao and Rams investigated the following sets
$$L(\alpha)=\left\{ x\in J_T : \ \lim_{n\to \infty}\frac{1}{n}\sum_{k=1}^n 1_{I_1}(T^kx)1_{I_1}(T^{2k}x)=\alpha\right\}\ \ (\alpha\in [0,1]).$$
By adapting the method of \cite{PS}, they obtained the Hausdorff spectrum of the above level sets $L(\alpha)$.

We point out that the methods used in \cite{PS} and \cite{LR} seem inconvenient to be generalised to other IFSs with many branches and more general potentials $\varphi$. Some more adaptive methods are needed to generalise Liao--Rams' results.
The aim of this paper is to use similar arguments as  in \cite{FSW1} to extend Liao--Rams' results to
the situation that we describe below.

Let $I_0,\cdots,I_{m-1}\subset[0,1]$ be $m$ intervals with disjoint interiors. Let $T: \cup_{i=0}^{m-1}I_i\to [0,1]$ be such that the restriction $T_{|I_i}$ is bijective and linear with slop $e^{\lambda_i}$, $ \lambda_i>0$ ($0\leq i\leq m-1$). Denote by $J_T$ the repeller of $T$.

Let $\ell\geq 2$ be an integer, and $\varphi$ be a function defined on $[0,1]^\ell$ taking real values. We assume that $\varphi$ is locally constant
in the sense that $\varphi$ is constant on each hyper-rectangle   $I_{i_1}\times I_{i_2}\times \cdots \times I_{i_\ell}$ $(0\leq i_1,i_2,\cdots,i_\ell \leq m-1)$.
With an abuse of notation, we write $$
\varphi(a_1,a_2,\cdots,a_\ell)=\varphi(i_i,i_2,\cdots,i_\ell)$$ for all $(a_1,a_2,\cdots,a_\ell)\in I_{i_1}\times I_{i_2}\times \cdots \times I_{i_\ell}$.

In this paper, we would like to study the following sets
$$L_\varphi(\alpha):=\left\{x\in J_T : \ \lim_{n\to \infty}\frac{1}{n}\sum \varphi(T^kx,T^{kq}x, \cdots, T^{kq^{\ell-1}}x)=\alpha\right\}, \alpha\in \R.$$ 
Our aim is to determine the Hausdorff dimension of $L_\varphi(\alpha)$. 

For simplicity of notations, we restrict ourselves to the case $\ell=2$ (the same arguments work for arbitrary $\ell\geq 2$ without any problem).
For any $s,r\in \R$, consider the non-linear transfer operator $\mathcal{N}_{(s,r)}$ on $\R_+^m$ defined by 
\begin{equation}\label{transfer equation}
\left(\mathcal{N}_{(s,r)}\underline{t} \right)_i=\left(\sum_{j=0}^{m-1}e^{s\varphi(i,j)-r\lambda_j}t_j\right)^{1/q},\ \ (i=0,\cdots,m-1).
\end{equation}
for all $\underline{t} =(t_j)_{j=0}^{m-1}\in\R_+^m$. 
In \cite{FSW1}, a family of similar  operators $\mathcal{N}_{s}\  (s\in\R)$ was defined. 
Notice that the Lyapunov exponents $\lambda_j$'s are now introduced
in the definition of  $\mathcal{N}_{(s,r)}$.  It will be shown in  Proposition \ref{Proposition  transfer} (see Section \ref{section Non-linear transfer})  that the equation $\mathcal{N}_{(s,r)}\underline{t} =\underline{t} $ admits a unique strictly positive solution $(t_0(s,r),\cdots,t_{m-1}(s,r))$. We then define the pressure function by
$$P(s,r)=(q-1)\log\sum_{j=0}^{m-1}t_j(s,r)e^{-r\lambda_j}.$$
 It will also be shown (Proposition \ref{Proposition  transfer}) that $P$ is real-analytic and convex, and even strictly convex if $\varphi$ is not constant and the $\lambda_j$'s are not all the same.

Let $A$ and $B$ be the infimum and the supremum respectively of the set $$\left\{a\in \mathbb{R} :\ \exists (s,r)\in \R^2 \ {\rm such\ that } \ \frac{\partial P}{\partial s}(s,r)=a\right\}.$$
Let $D_{\varphi}=\left\{\alpha\in \R : L_{\varphi}(\alpha)\neq \emptyset\right\}.$ Our main result is as follows.

\begin{theorem}\label{thm principal}
	Under the assumptions made above, we have 
\begin{enumerate}
\item[(i).] We have $D_{\varphi}=[A,B]$.  
\item[(ii).]  For any $\alpha\in (A,B)$, there exists a unique solution $(s(\alpha),r(\alpha))\in \R^2$ to the system
\begin{equation}\label{critical equation}
\left\{ \begin{array}{ll}
P(s,r)&=\alpha s\\
\frac{\partial P}{\partial s}(s,r)&=\alpha.
\end{array} \right.
\end{equation}
Furthermore, $\ s(\alpha)$ and $r(\alpha)$ are real-analytic functions of $\alpha \in (A,B)$.
\item[(iii).] The following limits exist: 
$$r(A):=\lim_{\alpha\downarrow A}r(\alpha), \ \ r(B):=\lim_{\alpha\uparrow B}r(\alpha).$$ 
\item[(iv).] For any $\alpha\in [A,B]$, we have
$$\dim_HL_\varphi(\alpha)=r(\alpha).$$

\end{enumerate}
\end{theorem}


The paper is organized as follows.
In Section 2, we first prove that the non-linear transfer operator $\mathcal{N}_{(s,r)}$ admits a unique positive fixed point $t(s,r)$ which is real-analytic and convex as a function of  $(s,r)$. Then we recall the class of {\em telescopic product measures} studied in \cite{KPS1, FSW1}. From each fixed point $t(s,r)$, we construct a special 
telescopic product measure, which will play the role of a Gibbs measure in our study of $L_\varphi(\alpha)$.
In Section 3, we study the local dimensions of the telescopic product measures defined by $t(s,r)$ and the formula of their local dimensions will be given.  Section 4 is devoted to the proof of (ii) of Theorem \ref{thm principal}. The assertions (i), (iii) and (iv) of  Theorem  \ref{thm principal} are proven in Section 5.

\section{Non-linear transfer equation and a class of special telescopic product measures. } \label{section Non-linear transfer}
Recall that $S=\{0,1,\cdots, m-1\}$ and $\Sigma_m=S^{\N^*}$. For $i\in S$,  let $f_i:[0,1]\to I_i$ be the branches of $T^{-1}$.  Define the coding map $\Pi: \Sigma_m\to [0,1]$ by
$$\Pi((x_k)_{k=1}^\infty)=\lim_{n\to\infty}f_{x_1}\circ f_{x_2}\cdots f_{x_n}(0).$$
Then we have $\Pi(\Sigma_m)=J_T$.
Define the subset $E_\varphi(\alpha)$ of $\Sigma_m$ which was studied in \cite{FSW,FSW1}:
\begin{equation*}
E_\varphi(\alpha):=\left\{(x_k)_{k=1}^\infty\in \Sigma_m : \ \lim_{n\to \infty}\frac{1}{n}\sum_{k=1}^n \varphi(x_k,x_{kq})=\alpha\right\}.
\end{equation*}
Then with a difference of a countable set, we have
$L_\varphi(\alpha)=\Pi(E_\varphi(\alpha)).$

In \cite{FSW,FSW1}, a family of Gibbs-type measures called  {\em telescopic product measures}  were used to compute the Hausdorff dimension of $E_\varphi(\alpha)$. Here  we construct a similar class of measures in order to determine the Hausdorff dimension of $L_\varphi(\alpha)$.
In the following, we suppose that $\varphi$ is not constant (otherwise the problem is trivial) and that the $\lambda_j$'s are not the same (otherwise the problem is reduced to the case considered in \cite{FSW,FSW1}).

\medskip

\subsection{Non-linear transfer operator}

In this subsection, we present some properties of the non-linear transfer operator $\mathcal{N}_{(s,r)}$, which will be used later.

\begin{proposition}\label{Proposition  transfer}
For any $s,r\in \R$,  the equation $\mathcal{N}_{(s,r)}y=y$ admits a unique solution $\underbar{{\rm t}}(s,r)=(t_0(s,r),\cdots,t_{m-1}(s,r))$ with strictly positive components, which  can be obtained as the limit of the iteration $\mathcal{N}_{(s,r)}^n\underbar{{\rm 1}}=: \underbar{\rm t}^n(s,r)$, where $\underbar{{\rm 1}}=(1,1,\cdots,1)$. The functions $t_i(s,r)$ and the pressure function $P(s,r)$ are real-analytic and strictly convex on $\R^2$.
\end{proposition}

\begin{proof}
(1). {\em Existence and uniqueness of solution.} Since $e^{s\varphi(i,j)+r\lambda_j} >0 $ for all $0\le i,j\le m-1$, the existence and uniqueness of solution are deduced directly from the following lemma.

\begin{lemma}{ \cite[Theorem 4.1]{FSW1} } \label{lemma existence and uniqueness}
For any matrix $A=(A(i,j))_{0\leq i,j\leq m-1}$ with strictly positive entries, there exists a unique fixed vector $\underbar{{\rm x}}=(x_0,\cdots,x_{m-1})\in \R^m$ with strictly positive components to the operator $\mathcal{N} : \R^m_+\to \R^m_+$ defined by 
$$\forall \underbar{\rm y}\in \R^m_+, \quad (\mathcal{N}\underbar{\rm y})_i=\left(\sum_{j=0}^{m-1}A(i,j)y_j\right)^{1/q}, \ (i=0,\cdots, m-1).$$
Furthermore, the fixed vector $\underbar{{\rm x}}$ can be obtained as $\underbar{{\rm x}}=\lim_n\mathcal{N}^n(\underbar{{\rm 1}})$.
\end{lemma}


(2). {\em Analyticity of $(s,r)\mapsto \underbar{{\rm t}}(s,r)$.} 
This has been proven in \cite[Proposition 4.2]{FSW1} for the case when all $\lambda_i$'s are the same. 
We adapt the proof given there with minor modifications.  

We consider the map $G: \R^2\times \R_+^{m}\to \R^{m}$ defined by 
$$\forall \underbar{z}=(z_0,\cdots,z_{m-1}) \quad G((s,r),\underbar{z})=\left(G_{i}((s,r),\underbar{z})\right)_{i=0}^{m-1},$$
where $$G_{i}((s,r),\underbar{z})=z_{i}^q-\sum_{j=0}^{m-1}e^{s\varphi(i,j)-r\lambda_j}z_{j}.$$
It is clear that $G$ is real-analytic. By Lemma \ref{lemma existence and uniqueness}, for any fixed $(s,r)\in \R^2$, $\underbar{t}(s,r)$ is the unique positive vector satisfying $$G((s,r),\underbar{t}(s,r))=0.$$
By the Implicit Function Theorem, to prove the analyticity of $(s,r)\mapsto \underbar{t}(s,r)$, we only need to show that the Jacobian matrix 
$$M(s)=\left(\frac{\partial G_{i}}{\partial z_j}\left((s,r),\underbar{t}(s,r)\right)\right)_{0\le i,j \le m-1}$$ is invertible  for all $(s,r)\in\R^2$. To this end, we consider the following matrix
$$\widetilde{M}(s)=\left(t_j(s,r) \frac{\partial G_{i}}{\partial z_j}((s,r),\underbar{t}(s,r)) \right)_{0\le i,j \le m-1},
$$ 
obtained by multiplying the $j$-th column of $M(s)$ by $t_j(s,r)$ for each $0\le j\le m-1$.
Then  $\det (M(s))\neq0$ if and only if  $\det(\widetilde{M}(s))\neq 0$.
So it suffices to prove that $\widetilde{M}(s)$ is invertible. We will show that $\widetilde{M}(s)$ is strictly diagonal dominating. Then by Gershgorin Circle Theorem (see e.g. \cite[Theorem 1.4, page 6]{Varga}), $\widetilde{M}(s)$ is invertible. 

Recall that  a matrix is said to be strictly diagonal dominating if for every row of the matrix, the modulus
 of the diagonal entry in the row is strictly larger than the sum of the modulus of all the other (non-diagonal) entries in that row.

Now we are left to show that for any $0\le i\le m-1$, 
\begin{equation}\label{analyticity 1}
\left|t_i(s,r) \frac{\partial G_{i}}{\partial z_i}((s,r),\underbar{t}(s,r))\right|- \sum_{\substack{0\le j\le m-1\\ i\neq j}}\left|t_j(s,r) \frac{\partial G_{i}}{\partial z_j}((s,r),\underbar{t}(s,r))\right|>0.
\end{equation}
In fact, we have
$$
\frac{\partial G_i}{\partial z_j}((s,r),\underbar{t}(s,r)) = \left\{ \begin{array}{ll}
qt_i^{q-1}(s,r)-e^{s\varphi(i,i)-r\lambda_i} & \textrm{ if } j=i,\\
e^{s\varphi(i,j)+r\lambda_j} & \textrm{otherwise.}
\end{array} \right.
$$
Then, substituting the last expression into (\ref{analyticity 1}), we deduce that the left hand side of (\ref{analyticity 1}) is equal to
\begin{equation}\label{equation analyticity 1}
qt_i^{q}(s,r)-\sum_{j=0}^{m-1}e^{s\varphi(i,j)-r\lambda_j}t_j(s,r).
\end{equation}
By the fact that $\underbar{t}(s,r)$ is the fixed vector of $\mathcal{N}_{(s,r)} $, \eqref{equation analyticity 1} is equal to $(q-1)t_i^q(s,r)$ which is strictly positive.

(3). {\em Convexity of $ \underbar{\rm t}(s,r)$ and $P(s,r)$.}
When all $\lambda_i$'s are the same, the convexity results of $ \underbar{t}(s,r)$ and $P(s,r)$ have been proven in detail in Sections 4 and 5 of \cite{FSW1} by studying the operator $\mathcal{N}_{(s,r)}$.  The main idea there is to prove by induction the convexity of each $(s,r)\mapsto \underbar{t}^n(s,r)$. Then the limit $\underbar{t}(s,r)=\lim_n \underbar{t}^n(s,r)$ is also convex. For the strictly convexity of $\underbar{t}(s,r)$, one uses analyticity property and the fact that a convex analytic function is either strictly convex or linear. 

We will omit the proofs which are elementary and are just minor modifications of those of \cite{FSW1}. One can refer to Sections 4, 5 and also 10 of \cite{FSW1}.
\end{proof}

To end this subsection, we give the following remark on the monotonicity of the function $r\mapsto P(s,r)$.

\begin{remark}\label{remark monotone of P}
Observe that for any fixed $s\in \R$, the function $r\mapsto \mathcal{N}_{(s,r)}^n\bar{1}$ is decreasing for all $n$. Thus for all $0\leq i \leq m-1$, the function $r\mapsto t_i(s,r)$ is also decreasing, and so is the function $r\mapsto P(s,r)$. 
\end{remark}

\medskip

\subsection{Construction of telescopic product measures and law of large numbers}

An important tool for the study of the multiple ergodic average of $\varphi$, introduced in \cite{KPS,KPS1} and used in \cite{FSW,FSW1,PS, LR},  is the telescopic product measure. This class of measures will also be  the  main ingredient of our proofs concerning the estimate of Hausdorff dimension of $L_\varphi(\alpha)$. 
Let us recall the definition of the telescopic product measure.  Consider the following partition of $\N^*$:
 $$
 \N^*=\bigsqcup_{i\geq 1,q\nmid i}\Lambda_i\ \ {\rm with}\ \Lambda_i=\{iq^j\}_{j\ge 0}.
 $$
 Then we decompose $\Sigma_m$ as follows:
 $$
 \Sigma_m=\prod_{i\geq1,q\nmid i}S^{\Lambda_i}.
 $$

 Let $\mu$ be a probability measure on $\Sigma_m$.
 We consider $\mu$ as a measure on $S^{\Lambda_i}$, which is identified with $\Sigma_m$, for every $i$ with $q\nmid i$.  Let
$\mu_i$ be a copy of $\mu$ on $S^{\Lambda_i}$ and $\P_{\mu}=\prod_{i\leq n,q\nmid i}\mu_i$.
 More precisely,  for any word $u$ of length $n$ we define
 $$
 \P_{\mu}([u])=\prod_{i\leq n,q\nmid i}\mu([u_{|_{\Lambda_i}}]),
 $$
 where $[u]$ denotes the cylinder set of all sequences starting with $u$ and 
 $$u_{|_{\Lambda_i}}=u_iu_{iq}\cdots u_{iq^j}, \ iq^j\leq |u|<iq^{j+1}.$$

Below, we construct
a special class of Markov measures whose initial laws and  transition probabilities
 are determined by the fixed point $(t_i(s,r))_{i\in S}$ of the operator $\mathcal{N}_{(s,r)}$.
The corresponding telescopic product measure will play a central role in the study of $E_\varphi(\alpha)$.

Recall that $(t_i(s,r))_{i\in S}$ satisfies 

$$
t_i(s,r)^q=\sum_{j=0}^{m-1}e^{s\varphi(i,j)-r\lambda_j}t_j(s,r),\ \ (i=0,\cdots,m-1).
$$
The functions $t_i(s,r)$ allow us to define a Markov measure
$\mu_{s,r}$ with initial law $\pi_{s,r}=(\pi(i))_{i\in S}$ and  probability transition matrix
$Q_{s,r}=(p_{i,j})_{S\times S}$  defined by
\begin{equation}\label{construction of special Markov measure}
        \pi(i) = \frac{t_i(s,r)e^{-r\lambda_i}}{t_0(s,r)e^{r\lambda_0}  + \cdots + t_{m-1}(s,r)e^{r\lambda_{m-1}}},
        \qquad p_{i, j}= e^{s \varphi(i, j)-r\lambda_j} \frac{t_j(s,r)}{t_i(s,r)^q}.
\end{equation}
We denote by $\mathbb{P}_{s,r}$ the telescopic product measure associated
to $\mu_{s,r}$.  Recall  that $\Pi$ is the coding map from $\Sigma_m$ to $[0,1]$.
Define $$\nu_{s,r}=\Pi_*\mathbb{P}_{s,r}=\mathbb{P}_{s,r}\circ\Pi^{-1}.$$

We will use the following law of large  numbers which is proved in \cite{FSW1}.

\begin{theorem}[Theorem 2.6 in \cite{FSW1}] \label{thm esperence general formula}
Let $\mu$ be any probability measure on $\Sigma_m$ and let
$F$ be a real-valued function defined  on $S\times S$.  For $\P_{\mu}$ a.e. $x\in \Sigma_m$
we have
$$\lim_{n\to\infty}\frac{1}{n}\sum_{k=1}^nF(x_k, x_{kq})=(q-1)^2\sum_{k=1}^\infty\frac{1}{q^{k+1}}\sum_{j=0}^{k-1}\E_\mu F(x_j,x_{j+1}).$$
\end{theorem}

\medskip

\section{Local dimension of $\nu_{s,r}$}
For a Borel measure $\mu$ on a metric space $X$, the lower local dimension of $\mu$ at a point $x\in X$ is defined by
\[
\underline{D}(\mu,x):=\liminf_{r\to 0} \frac{\log \mu(B(x,r))}{\log r}.
\]
If the limit exists, then the limit will be called the local dimension of $\mu$ at $x$, and denoted by ${D}(\mu,x)$.

In this section, we study the local dimension of $\nu_{s,r}$. The main results of this section are 
Propositions \ref{proposition local dimension level set}, \ref{proposition support on level set}, and \ref{proposition exact and formula}. Proposition \ref{proposition local dimension level set} gives estimates of the local dimensions of $\nu_{s,r}$ on the level set $L_\varphi(\alpha)$. Proposition \ref{proposition support on level set} proves that  $\nu_{s,r}$  is supported on $L_\varphi(\frac{\partial P}{\partial s}(s,r))$. In Proposition \ref{proposition exact and formula}, it is shown that $\nu_{s,r}$ is exact dimensional, i.e., the local dimension of $\nu_{s,r}$ exists and is constant almost surely. The exact formula of this constant is given as well.

We first give an explicit relation between the mass $\P_{s,r}([x_1^n])$ and the multiple ergodic sum $\sum_{j=1}^n\varphi(x_j,x_{qj})$. For $x\in \Sigma_m$, define 
$$B_n(x)=\sum_{j=1}^n \log t_{x_j}(s,r).$$
\begin{proposition}\label{prop basic formula loc dim} We have
$$
\log\P_{s,r}([x_1^n])=s\sum_{j=1}^{\lfloor
\frac{n}{q}\rfloor}\varphi(x_j,x_{jq})-(n-\lfloor \frac{n}{q}\rfloor)\frac{P(s,r)}{q-1}-r\sum_{j=1}^n\lambda_{x_j}-qB_{\lfloor
\frac{n}{q}\rfloor}(x)+B_n(x).
$$
\end{proposition}

\begin{proof}
For $q \nmid i$, let $\Lambda_i(n)=\Lambda_i\cap[1,n]$. By the definition of  $\P_{s,r}$, we have
\[-\log \P_{s,r}[x_1^n]=-\sum_{q\nmid i,i\leq n}\log\mu_{s,r}[x_1^n|_{\Lambda_i(n)}]. \]

We classify $\Lambda_i(n)$ ($q\nmid i, i\leq n$) according to their
length $|\Lambda_i(n)|$. 
We have $\min_{q\nmid i, i\leq n} |\Lambda_i(n)|=1$ and $\max_{q\nmid i, i\leq n} |\Lambda_i(n)|=\lfloor\log_q n\rfloor$.
 Observe that $|\Lambda_i(n)|=k$ if and only if $\frac{n}{q^k}<i\leq \frac{n}{q^{k-1}}.$ 
So
\begin{equation} \label{eq:section33}
-\log\P_{s,r}[x_1^n]=- \sum_{k=1}^{\lfloor \log_q
n\rfloor}\sum_{\frac{n}{q^k} <i\leq \frac{n}{q^{k-1}}, q\nmid
i}\log\mu_{s,r}[x_1^n|_{\Lambda_i(n)}].
\end{equation}
Denote $t_\emptyset(s,r):=\sum_{j\in S}t_j(s,r)e^{-r\lambda_j}$. For simplicity, we also write $t_\emptyset$ and $t_j$ for $t_\emptyset(s,r)$ and $t_j(s,r)$ and keep their dependences on $s$ and $r$ in mind.

By the definition of $\mu_{s,r}$, for $i$ with $\frac{n}{q^k} <i\leq \frac{n}{q^{k-1}}$, we have
\begin{align*} 
\log\mu_{s,r}[x_1^n|_{\Lambda_i(n)}] 
  &=  \log\frac{t_{x_i} e^{-r\lambda_{x_i}}}{t_\emptyset}+\sum_{j=1}^{k-1}\log
 \left(e^{s\varphi(x_{iq^{j-1}},x_{iq^j})-r\lambda_{x_{iq^j}}}\frac{t_{x_{iq^j}}}{t_{x_{iq^{j-1}}}^q}\right)\\
&=s S_{n,i}\varphi(x)-(q-1)S_{n,i}t(x)+\log
t_{x_{iq^{k-1}}}-rS_{n,i}\lambda(x)-\log t_\emptyset,
\end{align*}
where
$S_{n,i}\varphi(x)=\sum_{j=1}^{k-1}\varphi(x_{iq^{j-1}},x_{iq^j})$,
$S_{n,i}t(x)=\sum_{j=1}^{k-1}\log t_{x_{iq^{j-1}}}$ and $S_{n,i}\lambda(x)=\sum_{j\in\Lambda_i(n)} \lambda_{x_{j}}$. Substituting the above expressions in 
(\ref{eq:section33}) and noticing that $\frac{n}{q^k}<i\leq \frac{n}{q^{k-1}}$ is equivalent to $\frac{n}{q}<iq^{k-1}\leq n$, we obtain
\begin{eqnarray*}
\log\P_{s,r}[x_1^n] &=&s \sum_{q\nmid i,i\leq
n}S_{n,i}\varphi(x)-(q-1)\sum_{q\nmid i,i\leq
n}S_{n,i}t(x) +\sum_{\frac{n}{q} \leq \ell<n}\log t_{x_\ell}\\
&&-r \sum_{q\nmid i,i\leq
n}S_{n,i}\lambda(x)-\sharp\{q\nmid i,i\leq n\}\log
t_\emptyset\\
&=&s \sum_{j=1}^{\lfloor \frac{n}{q}\rfloor}\varphi(x_j,x_{jq})-(q-1)\sum_{j=1}^{\lfloor \frac{n}{q}\rfloor}
\log t_{x_j}+\sum_{\ell=\lfloor \frac{n}{q}\rfloor+1}^n\log t_{x_\ell}\\
&&-r\sum_{j=1}^n\lambda_{x_j}-(n-\lfloor
\frac{n}{q}\rfloor)\log t_\emptyset.
\end{eqnarray*}
We then end the proof by observing that $(q-1)\log t_\emptyset(s,r)=P(s,r)$ and 
$$-(q-1)\sum_{j=1}^{\lfloor \frac{n}{q}\rfloor}
\log t_{x_j}
+\sum_{\ell=\lfloor \frac{n}{q}\rfloor+1}^n\log t_{x_\ell}=-qB_{\lfloor
\frac{n}{q}\rfloor}(x)+B_n(x).$$
\end{proof}

\subsection{Local dimensions of $\nu_{s,r}$ on level sets.}
As an application of Proposition \ref{prop basic formula loc dim}, we obtain an upper bound for the local dimension of  $\nu_{s,r}$ on $L_\varphi(\alpha)$ in Proposition \ref{proposition local dimension level set} below. The following elementary result will be useful for the estimates of local dimension
of $\nu_{s,r}$.

\begin{lemma}\label{upper bounds elementary lemma}
Let $(a_n)_{n\geq 1}$ be a bounded sequence of non-negative real numbers. Then
\[\liminf_{n\to\infty}\left(a_{\lfloor \frac{n}{q} \rfloor}-a_n\right)\leq 0.\]
\end{lemma}

\begin{proof}
Let $b_l=a_{q^{l-1}}-a_{q^{l}}$ for
$l\in\N^*$. Then the boundedness implies
$$\lim_{l \to \infty }\frac{b_1+\cdots+b_l}{l}= \lim_{l \to \infty } \frac{a_1-a_{q^l}}{l}=0.$$
This in turn implies  $\liminf_{l\to \infty }b_l\leq 0$. Thus
$$\liminf_{l\to \infty }\left(a_{\lfloor \frac{n}{q} \rfloor}-a_n\right)\leq \liminf_{l\to \infty }b_l\leq 0.$$
\end{proof}

\begin{proposition}\label{proposition local dimension level set}
For any $x\in E_\varphi(\alpha)$, we have 
$$\liminf_n\frac{\log\nu_{s,r}(\Pi[x_1^n])}{\log|\Pi[x_1^n]|}\leq r+\limsup_n\frac{P(s,r)/q-\alpha s/q }{(\sum_{j=1}^n\lambda_{x_j})/n}.$$

\end{proposition}
\begin{proof}
Since  $\nu_{s,r}(\Pi[x_1^n])=\P_{s,r}([x_1^n])$, by Proposition \ref{prop basic formula loc dim} we can write $\log\nu_{s,r}(\Pi[x_1^n])$ as
$$
s\sum_{j=1}^{\lfloor
\frac{n}{q}\rfloor}\varphi(x_j,x_{jq})-(n-\lfloor \frac{n}{q}\rfloor)\frac{P(s,r)}{q-1}-r\sum_{j=1}^n\lambda_{x_j}-qB_{\lfloor
\frac{n}{q}\rfloor}(x)+B_n(x).$$
On the other hand, 
$\log|\Pi[x_1^n]|=-\sum_{j=1}^n\lambda_{x_j}$.
Thus, for $x\in E_\varphi(\alpha)$
\begin{align*}
&\liminf_n\frac{\log\nu_{s,r}(\Pi[x_1^n])}{\log|\Pi[x_1^n]|}\\
\leq & \limsup_n\frac{P(s,r)/q-\alpha s/q }{(\sum_{j=1}^n\lambda_{x_j})/n}+r+\liminf_n\frac{\frac{q}{n}B_{\lfloor \frac{n}{q}\rfloor}(x)-\frac{1}{n}B_n(x)}{(\sum_{j=1}^n\lambda_{x_j})/n}.
\end{align*}
Then, we end the proof by applying Lemma \ref{upper bounds elementary lemma} to the sequence $\frac{1}{n}B_n(x)$:
$$\liminf_n \left(\frac{q}{n}B_{\lfloor \frac{n}{q}\rfloor}(x)-\frac{1}{n}B_n(x)\right)\leq 0.$$
\end{proof}

\begin{remark} \label{remark  local dimension level set}
Denote $\lambda_{\min}=\min_i\lambda_i$ and $\lambda_{\max}=\max_i\lambda_i$. Let 
\[\widetilde{\lambda}(x):= \liminf_{n\to \infty} \frac1n\sum_{j=1}^n\lambda_{x_j}.\]
Then $\widetilde{\lambda}(x)\in [\lambda_{\min},\lambda_{\max}]$ and 
 \begin{align*}
 \limsup_n\frac{P(s,r)/q-\alpha s/q }{(\sum_{j=1}^n\lambda_{x_j})/n} &=\frac{P(s,r)/q-\alpha s/q }{\widetilde{\lambda}(x)}.
 \end{align*}
  So we deduce from Proposition \ref{proposition local dimension level set} that for any $x\in L_\varphi(\alpha)$
\begin{equation}\label{equation remark  local dimension level set 1}
\underline{D}(\nu_{s,r},x)\leq r+\frac{P(s,r)/q-\alpha s/q }{\widetilde{\lambda}(x)}, \ (s,r)\in \R^2.
\end{equation}
\end{remark}

We have estimated the local dimension of $\nu_{s,r}$ on the level set $L_\varphi(\alpha)$. In the following proposition we show that $\nu_{s,r}$ is supported on $L_\varphi(\frac{\partial P}{\partial s}(s,r))$.

\begin{proposition}\label{proposition support on level set}
For $\P_{s,r}$-a.e. $x=(x_i)_{i=1}^\infty\in \Sigma_m$,
we have
\begin{equation}\label{equation proposition support on level set 1}
\lim_{n\to\infty}\frac{1}{n}\sum_{k=1}^n\varphi(x_k,x_{kq})=\frac{\partial P}{\partial s}(s,r).
\end{equation}
In particular, $\nu_{s,r}\left(L_\varphi \Big(\frac{\partial P}{\partial s}(s,r)\Big)\right)=1$.
\end{proposition}

\begin{proof} We first prove the statement \eqref{equation proposition support on level set 1}.
By Theorem  \ref{thm esperence general formula}, we have for $\P_{s,r}$-a.e. $x\in \Sigma_m$
\begin{equation} \label{eq:section316}
\lim_{n\to\infty}\frac{1}{n}\sum_{k=1}^n\varphi(x_k,x_{kq})=(q-1)^2\sum_{k=1}^\infty\frac{1}{q^{k+1}}\sum_{h=0}^{k-1}\E_{\mu_{s,r}} \varphi(x_h,x_{h+1}).
\end{equation}
Thus we only need to prove that the right hand side of (\ref{eq:section316}) equals to $\frac{\partial P}{\partial s}(s,r)$.
Observe that $\E_{\mu_{s,r}} \varphi(x_h,x_{h+1})$ can be expressed as  
$$\pi
Q^h\widetilde{Q}(\underbar{1}),$$
with 
$$ \pi =\left( \frac{t_i(s,r)e^{-r\lambda_i}}{t_\emptyset(s,r)}\right)_{i\in S},
        \qquad Q= \left(e^{s \varphi(i, j)-r\lambda_j} \frac{t_j(s,r)}{t_i(s,r)^q}\right)_{(i,j)\in
S\times S}$$ and
$$\widetilde{Q}=\left(e^{s\varphi(i,j)-r\lambda_j}\varphi(i,j)\frac{t_j(s,r)}{t_i^q(s,r)}\right)_{(i,j)\in
S\times S}.$$
Recall that $(t_i(s,r))_i$  is the fixed point of $\mathcal{N}_{s,r}$:
\begin{equation} \label{eq:section317}
t_i^q(s,r)=\sum_{j=0}^{m-1}e^{s\varphi(i,j)-r\lambda_j}t_j(s,r),\quad (i,j)\in S\times
S. \end{equation} 
Taking the derivative with respect to $s$ of both sides of (\ref{eq:section317}), we get
\[qt_i^{q-1}(s)\frac{\partial t_i}{\partial s}(s,r)=\sum_{j=0}^{m-1}\left(e^{s\varphi(i,j)-r\lambda_j}\varphi(i,j)t_j(s,r)+ e^{s\varphi(i,j)}\frac{\partial t_j}{\partial s}(s,r)\right).\]
Dividing both sides of the above equation by
$t^q_i(s,r)$, we obtain
\begin{equation} \label{eq:section318}
\sum_{j=0}^{m-1}e^{s\varphi(i,j)-r\lambda_j}\varphi(i,j)\frac{t_j(s,r)}{t_i^q(s,r)}=q\frac{\frac{\partial t_i}{\partial s}(s,r)}{t_i(s,r)}-
\sum_{j=0}^{m-1}e^{s\varphi(i,j)} \frac{\frac{\partial t_j}{\partial s}(s,r)}{t_i^q(s,r)}.
\end{equation}
Let $w$ and  $v$ be two vectors defined by
\[w=q\left(\frac{\frac{\partial t_0}{\partial s}(s,r)}{t_0(s)}, \dots , \frac{\frac{\partial t_{m-1}}{\partial s}(s,r)}{t_{m-1}(s)}\right)^t\] and 
\[v=\left(\sum_{j=0}^{m-1}e^{s\varphi(0,j)} \frac{\frac{\partial t_j}{\partial s}(s,r)}{t_0^q(s)}, \dots , \sum_{j=0}^{m-1}e^{s\varphi(m-1,j)} \frac{\frac{\partial t_j}{\partial s}(s,r)}{t_{m-1}^q(s)}\right).\]
Then, by (\ref{eq:section318}), we have \[\widetilde{Q}(\underbar{1})=w-v.\]
Observe that $Qw=qv$, so
\begin{equation} \label{eq:section319}
\begin{split}
\sum_{h=0}^{k-1}\pi Q^h\widetilde{Q}(\underbar{1}) & = \sum_{h=0}^{k-1}\pi Q^h(w-v)  \\
&=\pi w+q\sum_{h=1}^{k-1}\pi Q^{h-1}v-\sum_{h=0}^{k-1}\pi Q^hv  \\
&= \pi w+ qS_{k-1}- S_k,
 \end{split}
 \end{equation}
where we denote $S_k=\sum_{h=0}^{k-1}\pi Q^{h}v$ for $k\ge 1$  and $ S_{0}=0$. Denote by $\alpha(s)$ the right hand side of (\ref{eq:section316}). Observe that $S_k/q^k\to 0 $ when $k\to \infty$. Substituting (\ref{eq:section319}) in
(\ref{eq:section316}), we obtain
\begin{align*}
\alpha(s)&= (q-1)^2\sum_{k=1}^\infty\frac{1}{q^{k+1}} \left(\pi w+
qS_{k-1}- S_k\right)\\
&=(q-1)^2\sum_{k=1}^\infty\frac{1}{q^{k+1}}\pi w \\& =  \frac{q-1}{q}\pi w =(q-1)\frac{\sum_{j=0}^{m-1}\frac{\partial t_j}{\partial s}(s,r)e^{-r\lambda_j}}{t_\emptyset(s,r)}= \frac{\partial P}{\partial s}(s,r).
\end{align*}
Now we show that $\nu_{s,r}\left(L_\varphi \Big(\frac{\partial P}{\partial s}(s,r)\Big)\right)=1$.
The formula \eqref{equation proposition support on level set 1} implies $$\mathbb{P}_{s,r}\big(E_\varphi(\frac{\partial P}{\partial s}(s,r))\big)=1.$$ Hence
 \begin{align*}
\nu_{s,r}\left(L_\varphi \Big(\frac{\partial P}{\partial s}(s,r)\Big)\right)&=\mathbb{P}_{s,r}\left(\Pi^{-1}\Big(L_\varphi\big(\frac{\partial P}{\partial s}(s,r)\big)\Big)\right)\\
&=\mathbb{P}_{s,r}\left(E_\varphi\Big(\frac{\partial P}{\partial s}(s,r)\Big)\right)=1.
 \end{align*}
\end{proof}


Let  $\lambda(s,r)$ be the expected limit with respect to $\P_{s,r}$ of the average of the Lyapunov exponents $\frac{1}{n}\sum_{k=1}^n\lambda_{\omega_k}$ with $\omega \in \Sigma_m$. By Theorem \ref{thm esperence general formula}, we have
$$\lambda(s,r)=(q-1)^2\sum_{k=1}^\infty\frac{1}{q^{k+1}}\sum_{j=0}^{k-1}\E_{\mu_{s,r}} \lambda_{\omega_j}.$$
As an application of Proposition \ref{proposition support on level set}, we show that the measure $\nu_{s,r}$ is exact dimensional and we have the following formula for its dimension.

\begin{proposition}\label{proposition exact and formula}
For $\nu_{s,r}$-a.e. $x$ we have
$$D(\nu_{s,r},x)=r+\frac{P(s,r)-s\frac{\partial P}{\partial s}(s,r)}{q\lambda(s,r)}.$$
\end{proposition}

\begin{proof}
We only need to show that for $\P_{s,r}$-a.e. $y\in\Sigma_m$
\begin{equation}\label{equation proposition exact and formula 1}
\lim_{n\to\infty} \frac{\log \P_{s,r}([y_1^n])}{\log |\Pi ([y_1^n])|}=r+\frac{P(s,r)-s\frac{\partial P}{\partial s}(s,r)}{q\lambda(s,r)}.
\end{equation}
Since $|\Pi ([y_1^n])|=e^{-\sum_{k=1}^n\lambda_{y_k}}$, from the discussion preceding Proposition \ref{proposition exact and formula}, we get for  $\P_{s,r}$-a.e. $y$
\begin{equation}\label{equation proposition exact and formula 2}
\lim_{n\to\infty} \frac{\log |\Pi ([y_1^n])|}{n}=-\lambda(s,r).
\end{equation}
On the other hand, by Theorem \ref{thm esperence general formula}, Proposition \ref{prop basic formula loc dim} and Proposition \ref{proposition support on level set}, we have for $\P_{s,r}$-a.e. $y$
\begin{equation}\label{equation proposition exact and formula 3}
\lim_{n\to\infty} \frac{\log \P_{s,r}([y_1^n])}{n}= \frac{s}{q}\frac{\partial P}{\partial s}(s,r)-\frac{1}{q}P(s,r)-r\lambda(s,r).
\end{equation}
Combining \eqref{equation proposition exact and formula 2} and \eqref{equation proposition exact and formula 3}, we get \eqref{equation proposition exact and formula 1}.

\end{proof}

\section{Further properties of the pressure function and study of the system (\ref{critical equation})}

The main result of this section is Proposition \ref{proposition study system} below on the  solution of the system  (\ref{critical equation}). 

We will use the following lemma concerning the range of the partial derivatives of $P(s,r)$. 
Recall that $(A,B)=\left\{ \frac{\partial P}{\partial s}(s,r) : (s,r)\in \R^2\right\}$.
\begin{lemma}\label{lemma range derivation pressure}
For any $r\in\R$, we have
$$\left\{ \frac{\partial P}{\partial s}(s,r) : s\in \R\right\}=(A,B)$$
\end{lemma}

\begin{proof}
Fix $r_0\in \R$. Since $s\mapsto P(s,r)$ is convex, It suffices to show that
$$\lim_{s\to +\infty}\frac{\partial P}{\partial s}(s,r_0) =B \ \textrm{ and  }\ \lim_{s\to -\infty}\frac{\partial P}{\partial s}(s,r_0) =A.$$
We only give the proof for the case when $s$ goes to $+\infty$. The case for $s$ tending to $-\infty$ is similar. The proof will be done by contradiction. Suppose that there exists  $\epsilon >0$ such that 
$$\frac{\partial P}{\partial s}(s,r_0) \leq B-\epsilon  \ \textrm{ for all  }\  s\in \R.$$ 
By the  Mean Value Theorem, for any $s>0$, we have 
\begin{equation}\label{equation lemma range derivation1}
P(s,r_0)-P(0,r_0)\leq s(B-\epsilon).
\end{equation}
By the definition of $B$, there exists $(s',r')\in\R^2$ such that $\frac{\partial P}{\partial s}(s',r')=B-\epsilon/2$. By Proposition \ref{proposition support on level set}, $\nu_{s',r'}(L_\varphi(B-\epsilon/2))=1$, so $L_\varphi(B-\epsilon/2)\neq \emptyset.$ Let $x\in L_\varphi(B-\epsilon/2)$. By Proposition \ref{proposition local dimension level set} and Remark \ref{remark  local dimension level set}, we have 
$$\underline{D}(\nu_{s,r_0},x)\leq r_0+\frac{P(s,r_0)/q-(B-\epsilon/2) s/q }{\widetilde{\lambda}(x)}.$$ Substituting (\ref{equation lemma range derivation1}) in the above inequality, we get
$$\underline{D}(\nu_{s,r_0},x)\leq r_0+\frac{P(0,r_0)/q-\epsilon s/2q }{\widetilde{\lambda}(x)}. $$
Since $\widetilde{\lambda}(x)\in [\lambda_{\min},\lambda_{\max}]$, the second term in the right hand side of the above inequality tends to $-\infty$ when $s\to +\infty$. So, for $s$ large enough we must have $\underline{D}(\nu_{s,r_0},x)<0$. But this is impossible since $\nu_{s,r_0}$ is a probability measure. Thus, we conclude that $\lim_{s\to +\infty}\frac{\partial P}{\partial s}(s,r_0) =B $. 

\end{proof}

\begin{proposition}\label{proposition study system}
For any $\alpha\in (A,B)$, there exists a unique solution
$(s(\alpha),r(\alpha))\in\R^2$ to the system
\begin{equation}\label{critical equation1}
\left\{ \begin{array}{ll}
P(s,r)&=\alpha s\\
\frac{\partial P}{\partial s}(s,r)&=\alpha,
\end{array} \right.
\end{equation}
Moreover the functions $s(\alpha), r(\alpha)$ are analytic on
$(A,B)$.
\end{proposition}

\begin{proof}
1). {\em Existence and uniqueness of the solution $(s(\alpha),r(\alpha))$.} \  Fix $\alpha\in (A,B)$.  By Lemma \ref{lemma range derivation pressure} and the strict convexity of $s\mapsto P(s,r)$, for any $r\in \R$, there exists a unique $s=s(\alpha,r)\in \R$ such that
\begin{equation}\label{equation proposition study system 1}
\frac{\partial P}{\partial s}(s(\alpha,r),r)=\alpha.
\end{equation}
In the following, we will show that there exists a unique solution $r=r(\alpha)\in \R$ to the equation
$$P(s(\alpha,r),r)=\alpha s(\alpha,r).$$
Set $h(r):=P(s(\alpha,r),r)-\alpha s(\alpha,r)$. By \eqref{equation proposition study system 1}
\begin{eqnarray*} 
h'(r)
  &= & \frac{\partial P}{\partial s}(s(\alpha,r),r)\frac{\partial s(\alpha,r)}{\partial r}+\frac{\partial P}{\partial r}(s(\alpha,r),r)-\alpha \frac{\partial s(\alpha,r)}{\partial r}\\
&=& \frac{\partial P}{\partial r}(s(\alpha,r),r).
\end{eqnarray*}
For fixed $s$ the function $r\mapsto P(s,r)$ is strictly decreasing, since it is strictly convex and decreasing (Remark \ref{remark monotone of P}). So $ \frac{\partial P}{\partial r}(s(\alpha,r),r) <0$ and thus $h(r)$ is also strictly decreasing. For the rest of the proof, we only need to show  $\lim_{r\to+\infty} h(r)<0$ and $\lim_{r\to-\infty} h(r)>0$, then we conclude by applying the  Intermediate Value Theorem.

By Proposition \ref{proposition exact and formula}, we have 
$$\dim \nu_{s(\alpha,r),r} =r+\frac{P(s(\alpha,r),r)-s(\alpha,r)\alpha}{q\lambda(s(\alpha,r),r)}.$$
Observe that for any $r\in\R$, we have always $0\leq \dim \nu_{s(\alpha,r),r}\leq 1$ and $0<\lambda_{\min}\leq \lambda(s(\alpha,r),r)\leq \lambda_{\max}$. So we have
$$\lim_{r\to +\infty}h(r)=\lim_{r\to +\infty}\left(\dim \nu_{s(\alpha,r),r} -r\right)q\lambda(s(\alpha,r),r)<0.$$ 
Similarly,
$$\lim_{r\to -\infty}h(r)>0.$$

2). {\em  Analyticity of $(s(\alpha),r(\alpha))$.} \  Consider the map
$$F =
\left( \begin{array}{ccc}
F_1 \\
F_2
\end{array} \right)=
\left( \begin{array}{ccc}
P(s,r)- \alpha s\\
\frac{\partial P}{\partial s}(s,r)-\alpha
\end{array} \right).
$$
The jacobian matrix of $F$ is equal to
$$
J(F) :=
\left( \begin{array}{ccc}
\frac{\partial F_1}{\partial s} & \frac{\partial F_1}{\partial r} \\
\frac{\partial F_2}{\partial s} & \frac{\partial F_2}{\partial r} 
\end{array} \right)=
\left( \begin{array}{ccc}
\frac{\partial P}{\partial s}-\alpha & \frac{\partial P}{\partial r}\\
\frac{\partial^2 P}{\partial s^2} & \frac{\partial^2 P}{\partial r \partial s } 
\end{array} \right).
$$
So we have
$$\det (J(F))|_{s=s(\alpha),r=r(\alpha)}=-\frac{\partial^2 P}{\partial s^2}\cdot \frac{\partial P}{\partial r}\neq 0.$$
Then by the Implicit Function Theorem, $s(\alpha)$ and $r(\alpha)$ are analytic.
\end{proof}

\medskip

\section{Proof of Theorem \ref{thm principal}}

\subsection{Computation of $\dim_H L_\varphi(\alpha)$ for $\alpha\in (A,B)$}

We will use the following Billingsley Lemma.

\begin{lemma}[see e.g. Proposition 4.9. in \cite{Fal90}]\label{Billingsley}
Let $E\subset\Sigma_m$ be a Borel set  and let $\mu$ be a finite Borel measure on $\Sigma_m$.\\
\indent  {\rm (i)}  If $\mu(E) > 0$ and $\underline{D}(\mu,x)\geq d$ for $\mu$-a.e $x$, then  $\dim_H(E)\geq d$;\\
\indent {\rm (ii)} If $\underline{D}(\mu,x)\leq d$ for all $x \in E$, then   $\dim_H(E) \leq d$.
\end{lemma}

\medskip

\begin{theorem}\label{theorem dim of level set} 
For any $\alpha \in (A,B)$, we have 
$$\dim_HL_\varphi(\alpha)=r(\alpha).$$
\end{theorem}

\begin{proof}
By \eqref{equation remark  local dimension level set 1} and the equality  $P(s(\alpha),r(\alpha))=\alpha s(\alpha)$, we have
$$\underline{D}(\nu_{s(\alpha),r(\alpha)},x)\leq r(\alpha) \ \textrm{ for all }\ x\in L_\varphi(\alpha).$$
Then Lemma \ref{Billingsley} implies that 
$$\dim_HL_\varphi(\alpha)\leq r(\alpha). $$
By Proposition \ref{proposition support on level set}  and the equality  $\frac{\partial P}{\partial s}(s(\alpha),r(\alpha))=\alpha s(\alpha)$, we know that $$\nu_{s(\alpha),r(\alpha)}(L_\varphi(\alpha))=1.$$ 
On the other hand, by Proposition \ref{proposition exact and formula}, 
$$D(\nu_{s(\alpha),r(\alpha)},x)= r(\alpha) \ \textrm{ for }\ \nu_{s(\alpha),r(\alpha)}{\rm -a.e.} \  x.$$
Applying Lemma \ref{Billingsley} again, we obtain 
$$\dim_HL_\varphi(\alpha)\geq r(\alpha). $$
\end{proof}

\subsection{Range of $\{\alpha : L_\varphi(\alpha)\neq \emptyset\}$}

\begin{proposition} \label{range}
We have
$\{\alpha : L_\varphi(\alpha)\neq \emptyset\}\subset [A,B].$
\end{proposition}

\begin{proof} We prove it by contradiction.
Suppose that $L_\varphi(\alpha)\neq \emptyset$ for some
$\alpha>B$. Let $x\in
L_\varphi(\alpha)$. Then by  \eqref{equation remark  local dimension level set 1} and taking $r=0$, we
have
\begin{equation}\label{prop borne inf de mesure1}
\underline{D}(\nu_{s,0},x)\leq\frac{P(s,0)-\alpha s }{q\widetilde{\lambda}(x)} \ \textrm{ for all } s\in \R.
\end{equation}
On the other hand, by the mean value theorem, we have
\begin{equation}\label{prop borne inferieure de mes2}
P(s,0)-\alpha s=\frac{\partial P}{\partial s}(\eta_s,0)s-\alpha s+P(0,0)
\end{equation} for some real number $\eta_s$ between $0$ and $s$. In the following, we suppose that $s>0$.
Substituting (\ref{prop borne inferieure de mes2}) in (\ref{prop borne inf de mesure1}), we get 
$$\underline{D}(\nu_{s,0},x)\leq\frac{\frac{\partial P}{\partial s}(\eta_s,0)s-\alpha s+P(0,0)}{q\widetilde{\lambda}(x)}\leq \frac{(B -\alpha) s+P(0,0)}{q\widetilde{\lambda}(x)}.$$
Since $B-\alpha<0$ and $\widetilde{\lambda}(x)\in [\lambda_{\min},\lambda_{\max}]$, the last term in the above inequalities tend to $-\infty$ when $s\to +\infty$. But this is impossible since we have always $\underline{D}(\nu_{s,0},x)\geq 0$. Thus we must have $L_\varphi(\alpha)= \emptyset$ for any
$\alpha>B$. Similarly we can also prove that $L_\varphi(\alpha)= \emptyset$ for any
$\alpha<A$.
\end{proof}

As we will show, we actually have the equality $\{\alpha : L_\varphi(\alpha)\neq \emptyset\}= [A,B]$ (see Theorem \ref{theorem boundary pints}).

\medskip

\subsection{Computation of $\dim_H L_\varphi(A)$ and $\dim_H L_\varphi(B)$}

Now, we consider the level set $L_\varphi(\alpha)$ when $\alpha=A$ or $B$. The aim of this subsection is to prove the following theorem.
\begin{theorem}\label{theorem boundary pints}
{\rm(i).} The following limits exist: 
$$r(A):=\lim_{\alpha\to A}r(\alpha),\ \ r(B):=\lim_{\alpha\to B}r(\alpha).$$

{\rm(ii).} If $\alpha=A$ or $B$, then $L_\varphi(\alpha)\neq \emptyset$ and 
$$\dim_H L_\varphi(\alpha)=r(\alpha).$$

\end{theorem}

We will give the proof of Theorem \ref{theorem boundary pints} for the case $\alpha=A$, the proof for $\alpha=B$ is similar.

\subsubsection{Accumulation points of $\mu_{s(\alpha),r(\alpha)}$ when $\alpha $ tends to $A$.}
As all components of the vector $\pi_{s,r}$ and the matrix $Q_{s,r}$ (see formula \eqref{construction of special Markov measure}) are non-negative and bounded by 1, the set $\{(\pi_{s(\alpha),r(\alpha)},Q_{s(\alpha),r(\alpha)}), \alpha\in (A,B)\}$ is precompact. So there exists a sequence $(\alpha_n)_n\in (A,B) $ with $\lim_n\alpha_n=A$ such that the limits 
$$\lim_{n\to\infty}\pi_{s(\alpha_n),r(\alpha_n)}, \ \ \lim_{n\to\infty} Q_{s(\alpha_n),r(\alpha_n)}$$ exist. Using these limits as initial law and transition probability, we construct a  Markov measure which we denote by $\mu_\infty$. It
is clear that the Markov measure $\mu_{s(\alpha_n),r(\alpha_n)}$ corresponding to $\pi_{s(\alpha_n),r(\alpha_n)}$ and $Q_{s(\alpha_n),r(\alpha_n)}$
converges to $\mu_\infty$ with respect to the weak-star topology.
We denote by $\P_\infty$
 the telescopic product measure associated to $\mu_\infty$ and set $\nu_\infty:=\P_\infty\circ\Pi^{-1}$.

\begin{proposition} \label{Proposition boundary point full measure}
We have
$$\nu_\infty(L_\varphi(A))=1.$$
In particular, $L_\varphi(A)\neq \emptyset$.
\end{proposition}

\begin{proof}

Since $\nu_{\infty}(L_\varphi(A))=\P_\infty(E_\varphi(A))$, we only need to show that $\P_\infty(E_\varphi(A))=1$, i.e.,
 for $\P_{\infty}$-a.e. $x\in \Sigma_m$ we have
$$\lim_{n\to\infty}\frac{1}{n}\sum_{k=1}^n\varphi(x_k,x_{kq})=A.$$
By Theorem \ref{thm esperence general formula},  for $\P_{\infty}$-a.e. $x\in \Sigma_m$ the limit in the left hand side of the above equation equals $M(\mu_\infty)$ where $M$ is the functional on the space of probability measures defined by
$$M(\nu)=(q-1)^2\sum_{k=1}^\infty\frac{1}{q^{k+1}}\sum_{j=0}^{k-1}\E_\nu \varphi(x_j,x_{j+1}).$$
The function $\nu\mapsto M(\nu)$ is continuous, since the above series converges uniformly on $\nu$ and the function $\nu\mapsto \E_\nu \varphi(x_j,x_{j+1})$ is continuous for all $j$.
Since $\mu_{s(\alpha_n),r(\alpha_n)}$ converges to $\mu_{\infty}$ when $n\to \infty$, we have that$$\lim_{n\to \infty}M(\mu_{s(\alpha_n),r(\alpha_n)})=M(\mu_{\infty}).$$
Recall that the vector $(s(\alpha),r(\alpha))$ satisfies $\frac{\partial P}{\partial s}(s(\alpha),r(\alpha))=\alpha$.
By Proposition \ref{proposition support on level set}, we know that $$M(\mu_{s(\alpha_n),r(\alpha_n)})=\alpha_n.$$
So $$M(\mu_{\infty})=\lim_{n\to \infty}\alpha_n=A.$$

\end{proof}

From Theorem \ref{theorem dim of level set}, we know that  for each $\alpha\in (A,B)$, $r(\alpha)=\dim_HL_\varphi(\alpha)\in  [0,1]$. So, in particular the set $\{r(\alpha): \alpha\in (A,B)\}$ is bounded. 

We have the following formula for $\dim_H \nu_{\infty}$.

\begin{proposition}\label{Proposition boundary point dim formula}
The limit $r(A):=\lim_nr(\alpha_n)$ exists and 
we have 
$$\dim \nu_\infty =r(A).$$
\end{proposition}
\begin{proof}
Let $(\alpha_{n_k})_k$ be any subsequence of $(\alpha_n)_n$ such that the 
the limit $\lim_kr(\alpha_{n_k})$ exists. We will show that this limit is equal to  $\dim \nu_\infty$.

The measure $\nu_\infty$ is exact dimensional and its dimension is given by 
$$\dim \nu_\infty=\frac{\dim(\P_\infty)}{\lambda(\P_\infty)},$$
where $\dim(\P_\infty)$ is the a.e. local dimension of $\P_\infty$ and $\lambda(\P_\infty)$ is  the expected limit with respect to $\P_{\infty}$  of the average of the Lyapunov exponents $\frac{1}{n}\sum_{k=1}^n\lambda_{\omega_k}$ with $\omega \in \Sigma_m$, i.e.,
$$\lambda(\P_{\infty})=(q-1)^2\sum_{k=1}^\infty\frac{1}{q^{k+1}}\sum_{j=0}^{k-1}\E_{\mu_{\infty}} \lambda_{\omega_j}.$$
By similar arguments as used in the proof of Proposition \ref{Proposition boundary point full measure}, we can show that the functions 
$$\mu\mapsto  \dim(\P_\mu), \ \ \mu\mapsto  \lambda(\P_\mu)$$
are continuous on the space of probability measures.
Thus, we deduce that 
$$\dim \nu_\infty =\lim_{k\to\infty}\frac{\dim(\P_{\mu_{s(\alpha_{n_k}),r(\alpha_{n_k})}})}{\lambda(\P_{\mu_{s(\alpha_{n_k}),r(\alpha_{n_k})}})}=\lim_{k\to\infty}\dim \nu_{s(\alpha_{n_k}),r(\alpha_{n_k})} =\lim_{k\to\infty}r(\alpha_{n_k}),$$
where we have used Theorem \ref{theorem dim of level set} for the last equality.
Since the subsequence $(\alpha_{n_k})_k$ is arbitrary, we deduce that the limit $r(A):=\lim_nr(\alpha_n)$ exists and $\dim \nu_\infty =r(A)$.
\end{proof}

In the proof of Theorem \ref{theorem boundary pints}, we will use the following lemma. Recall that for $\alpha\in (A,B)$, the vector $(s(\alpha),r(\alpha))$ is the unique solution of the equation \eqref{critical equation1}.
\begin{lemma} \label{lemma boundary points}
There exists $A'\in (A,B)$ such that 
$$s(\alpha)<0 \ \ {\rm for }\ \alpha \in (A, A').$$
\end{lemma}
\begin{proof}
Let $$D:=\left\{\frac{\partial P}{\partial s}(0,r): r\in [0,1]\right\}.$$ Then $D$ is a compact set of $\R$. Since for any $r\in \R$ the function $s\mapsto \frac{\partial P}{\partial s}(s,r)$ is strictly increasing and  $\inf_{s\in \R}\frac{\partial P}{\partial s}(s,r)=A$ (Lemma \ref{lemma range derivation pressure}), we get
$\frac{\partial P}{\partial s}(0,r)>A$ for all $r\in \R$. Thus we have $A':=\min \{D\}>A$. Now, we consider the following subset of $D$:  
$$D':=\left\{\frac{\partial P}{\partial s}(0,r(\alpha)): \alpha\in (A,B)\right\}.$$
We have $\inf D'\geq A'>A$. 
For any $\alpha<A'$, we have 
$$\frac{\partial P}{\partial s}(s(\alpha),r(\alpha))=\alpha<A'\le \frac{\partial P}{\partial s}(0,r(\alpha)).$$
Using again the fact that  the function $s\mapsto \frac{\partial P}{\partial s}(s,r)$ is strictly increasing, we get
$$s(\alpha)<0 \ \ {\rm for}\ \ \alpha\in (A,A').$$
\end{proof}

Now, we can give the proof of Theorem \ref{theorem boundary pints}.
\begin{proof}[Proof of Theorem \ref{theorem boundary pints}]
(1). Fix any sequence $(\beta_n)_n\in (A,B)$ with $\lim_n\beta_n=A$. Then there exists a subsequence $(\beta_{n_k})_k $ of $(\beta_n)_n$  such that the limits 
$$\lim_{k\to\infty}\pi_{s(\beta_{n_k}),r(\beta_{n_k})}, \ \ \lim_{k\to\infty} Q_{s(\beta_{n_k}),r(\beta_{n_k})}$$ exist. With a same proof of Proposition \ref{Proposition boundary point dim formula}, we can show that the limit $\lim_kr(\beta_{n_k})$ exists and equals to $\dim \nu_\infty$.
Thus, we deduce that the limit $\lim_{\alpha\to A}r(\alpha)$ exists and equals to $\dim \nu_\infty$.

(2). We will show that $$\dim_H L_\varphi(A)=r(A).$$ 

By Proposition \ref{Proposition boundary point full measure} and \ref{Proposition boundary point dim formula} and 
Lemma \ref{Billingsley}, we get 
$$\dim_H L_\varphi(A)\geq r(A).$$
We now show the reverse inequality. By \eqref{equation remark  local dimension level set 1} and again Lemma \ref{Billingsley}, we obtain $$\dim_HL_\varphi(A)\leq r+\frac{P(s,r)-A s }{q\widetilde{\lambda}(x)} \ \ {\rm for\ any}\ (s,r)\in \R^2.$$
Note that $\widetilde{\lambda}(x)\in [\lambda_{\min},\lambda_{\max}]$, so in particular $\widetilde{\lambda}(x)>0$.
For any $\alpha\in (A,A')$, we have
$$P(s(\alpha),r(\alpha))-As(\alpha)=P(s(\alpha),r(\alpha))-\alpha s(\alpha)+(\alpha-A)s(\alpha)=(\alpha-A)s(\alpha)<0,$$
where for the second equality we have used the fact that $P(s(\alpha),r(\alpha))=\alpha s(\alpha)$ and  the last inequality follows from Lemma \ref{lemma boundary points}. Thus, we deduce that 
$$\dim_HL_\varphi(A)\leq r(\alpha) \ \ {\rm for\ all}\ \ \alpha\in (A,A']. $$
Since $\alpha_n\to A$ and $r(\alpha_n)\to r(A)$, we have 
$$\dim_HL_\varphi(A)\leq \lim_{n\to\infty}r(\alpha_n)=r(A).$$

\end{proof}

\bigskip

{\em Acknowledgement:} The first author is  partially supported by NSFC No. 11471132 and by the self-determined research funds of CCNU (No. CCNU14Z01002) from the basic research and operation of MOE. The third author acknowledges the support of Academy of Finland, the Centre of Excellence in Analysis and Dynamics Research.
\medskip

\bigskip
\bigskip


\end{document}